\documentclass{svmult}

\usepackage{makeidx}         
\usepackage{graphicx}        
\usepackage{multicol}        
\usepackage[bottom]{footmisc}
\usepackage{natbib}
\usepackage{url}
\usepackage[section]{algorithm}
\usepackage{algpseudocode}

\renewcommand{\PackageWarningNoLine}[2]{}
\usepackage{amsmath}
\usepackage{amsfonts}

\usepackage{booktabs} 
\usepackage{subfig}

\usepackage[latin9]{inputenc}

\usepackage[usenames,dvipsnames]{color}
\usepackage{soul}

\usepackage[final]{changes}

\definechangesauthor[name=Marcella,color=blue]{MB}
\definechangesauthor[name=report1,color=red]{R1}
\definechangesauthor[name=report1,color=blue]{R2}

\begin{document}

\title*{Two-level preconditioners for the Helmholtz equation}

\author{
  Marcella Bonazzoli\inst{1}\and
  Victorita Dolean\inst{1,2}\and
  Ivan G. Graham \inst{3} \and 
  Euan A. Spence \inst{3} \and
  Pierre-Henri Tournier \inst{4}}
  
\authorrunning{M. Bonazzoli, V. Dolean, I.G. Graham, E.A. Spence, P.-H. Tournier}  

\institute{
  \inst{1}
  Universit\'{e} C\^{o}te d'Azur, CNRS, LJAD, France,
  \email{marcella.bonazzoli@unice.fr}\\ 
  \inst{2}
  University of Strathclyde, Glasgow, UK,
  \email{Victorita.Dolean@strath.ac.uk}\\ 
  \inst{3} 
  University of Bath, UK,
  \email{I.G.Graham@bath.ac.uk, E.A.Spence@bath.ac.uk}\\
  \inst{4}
  UPMC Univ Paris 06, LJLL, Paris, France,
  \email{tournier@ljll.upmc.fr}
}

%
%
\maketitle

\abstract*{
In this paper we compare numerically two different coarse space definitions for two-level domain decomposition preconditioners for the Helmholtz equation, both in two and three dimensions. While we solve the pure Helmholtz problem without absorption, the preconditioners are built from problems with absorption. In the first method, the coarse space is based on the discretization of the problem with absorption on a coarse mesh, with diameter constrained by the wavenumber. In the second method, the coarse space is built by solving local eigenproblems involving the Dirichlet-to-Neumann (DtN) operator.
} 


\section{Introduction}
\label{bonazzoli_mini_10_sec:intro}

Solving the Helmholtz equation  $-\Delta u -k^2 u = f$ is a challenging task because of its indefinite nature and its highly oscillatory solution when the wavenumber $k$ is high. Although there have been different attempts to solve it efficiently, we believe that there is no established and robust preconditioner, whose behavior is independent of $k$, for general decompositions into subdomains. In \cite{Conen:2014:CSH} a two-level preconditioner was introduced, where the coarse correction involves local eigenproblems of Dirichlet-to-Neumann (DtN) maps. 
This method proved to be very robust with respect to heterogeneous coefficients compared to the reference preconditioner based on plane waves, and its construction is completely automatic without the need for parameter tuning. 
Another method was developed in \cite{graham:2015:domain,Graham:2016:RRD}, 
where two-level domain decomposition approximations of the Helmholtz equation with absorption $-\Delta u - (k^2+\mathtt{i}\varepsilon) u = f $ were used as preconditioners for the pure Helmholtz equation without absorption; there the coarse correction is based on a coarse mesh with diameter constrained by $k$.
Our purpose is to compare numerically the performance of the latter with the two-level method based on DtN maps, both in two and three dimensions. 

\section{Definition of the problem}
\label{bonazzoli_mini_10_sec:helmholtzEq}

Consider the interior Helmholtz problem of the following form: let $\Omega \subset \mathbb R^d$, $d=2,3$, be a polyhedral, bounded domain; find
$u \colon \Omega \to \mathbb C$ such that
\begin{subequations}
\label{bonazzoli_mini_10_eq:helmholtzEquation}
\begin{align}
	- \Delta u - (k^2+\mathtt{i}\varepsilon) u &= f && \mbox{in } \Omega, \\
	\frac{\partial u}{\partial n} - \mathtt{i} \eta u &= 0 && \mbox{on } \Gamma=\partial \Omega.
	\label{bonazzoli_mini_10_eq:helmholtzRobin}
\end{align} 
\end{subequations}
The wavenumber $k$ is given by $k(\mathbf{x}) = \omega/{c(\mathbf{x})}$, where $\omega$ is the angular frequency and $c$ is the speed of propagation that might depend on $\mathbf{x}\in \Omega$; we take $\eta=\mathrm{sign}(\varepsilon)k$ \added[id=R2]{if $\varepsilon \ne 0$,  $\eta=k$ if $\varepsilon=0$,} as Robin boundary condition parameter. We are interested in solving the problem in the case $\varepsilon=0$, using $\varepsilon$ as a parameter when building the preconditioner.
The variational formulation of Problem \eqref{bonazzoli_mini_10_eq:helmholtzEquation} is: find $u \in V = H^1(\Omega)$ such that 
$ 
	a_{\varepsilon}(u,v) = F(v), \, \forall v \in V,
$ 
where $a_{\varepsilon}(.,.)\colon V \times V \to \mathbb{C}$ and $F\colon V\to \mathbb{C}$ are defined
by
\begin{equation*}
	a_{\varepsilon}(u,v) = \int_{\Omega}{\left(\nabla u \cdot \overline{\nabla v} - (k^2+\mathtt{i}\varepsilon) u \overline{v}\right)}
	-
	\int_{\Gamma}{\mathtt{i} \eta u \overline{v} }, \quad
	F(v) = \int_{\Omega}{f \overline{v}}.
\end{equation*}
Note that if $\varepsilon \ne 0$ and $\eta=\mathrm{sign}(\varepsilon)k$, $a_{\varepsilon}$ is coercive (see \S 2 in \cite{graham:2015:domain}).
We consider a discretization of the variational problem
using piecewise linear finite elements on a uniform
simplicial mesh $\mathcal T_h$ of $\Omega$.
Denoting by $V_h \subset V$ the corresponding finite element space and by $\{\phi_k\}_{k=1}^n$ its basis functions,  $n := \text{dim}(V_h)$, the discretized problem reads:
find $u_h \in V_h$ such that 
$ a_{\varepsilon}(u_h,v_h) = F(v_h), \, \forall v_h \in V_h$,  
that is, in matrix form,
\begin{equation}
	\label{bonazzoli_mini_10_eq:matrixEqHelm}
	 A_{\varepsilon} \vec u = \vec f,
\end{equation}
where the coefficients of the matrix $ A_{\varepsilon} \in \mathbb{C}^{n
\times n}$ and the right-hand side $\vec f \in \mathbb{C}^n$ are given by 
$ (A_{\varepsilon})_{k,l} = a(\phi_l, \phi_k)$ and $(\vec f)_k = F(\phi_k)$.
The matrix $A_{\varepsilon} $ is complex, symmetric (but not Hermitian), and indefinite if $\varepsilon=0$.

\section{Two-level domain decomposition preconditioners}

In the following we will define the domain decomposition preconditioners for the linear system  $A_{0} \vec u = \vec f$ resulting from the discretization of the Helmholtz problem without absorption ($\varepsilon=0$). These are two-level Optimized Restricted Additive Schwarz (ORAS) algorithms, where ``optimized" refers to the use of Robin boundary conditions at the interface between subdomains.
In the terminology of \cite{graham:2015:domain}, the prefix O is replaced with Imp, which stands for impedance (i.e.~Robin) boundary conditions.

First of all, consider a decomposition of the domain $\Omega$ into a set of overlapping subdomains $\{\Omega_j\}_{j=1}^{N_\text{sub}}$, with each subdomain consisting of a union of elements of the mesh $\mathcal T_h$.  
Let $
	V_h(\Omega_j) = \left\{ v|_{\Omega_j}: v \in V_h \right\}$, $1\leq j \leq
	N_\text{sub},
$ 
denote the space of functions in $V_h$ restricted to the subdomain $\Omega_j$.
 \replaced[id=R1]{Let $n_j$ be the dimension of $V_h(\Omega_j)$, $1\leq j \leq N_\text{sub}$.}{Let $n_j := \# \, {\mathrm{dof}\left(\Omega_j\right)}$, $1\leq j \leq N_\text{sub}$, where 
$\mathrm{dof}(D) := \left\{ k: \mathrm{supp}\left(\phi_k\right) \subset \overline{D} \right\}$
represents the degrees of freedom (dofs) of $V_h(\Omega_j)$.}
For $1 \leq j \leq N_\text{sub}$, 
we define a restriction operator $\mathcal R_j\colon V_h \to V_h(\Omega_j)$ by
injection, i.e. for $u \in V_h$ we set
$\left(\mathcal R_j u\right)\left(\vec x_i\right) = u(\vec x_i)$ for all
$\vec x_i \in \Omega_j$. 
We denote  by $R_j$ the corresponding Boolean matrix in $\mathbb R^{n_j
\times n}$ that maps coefficient vectors of functions in $V_h$ to coefficient
vectors of functions in $V_h(\Omega_j)$.
Let $D_j \in
\mathbb{R}^{n_j \times n_j}$ be a diagonal matrix
corresponding to a partition of unity in the sense that
$
	\sum_{i=1}^{N_\text{sub}} \tilde R_i^T  R_i = I,
$
where  $\tilde R_j := D_j R_j$. Then the \emph{one-level ORAS} preconditioner (which is also the {one-level \replaced[id=R2]{ImpRAS}{ ImpHRAS}} of \cite{graham:2015:domain}) reads 
\begin{equation}
	\label{bonazzoli_mini_10_eq:ORAS}
	M^{-1}_{1,\varepsilon} := \sum_{j=1}^{N_\text{sub}} \tilde R_j^T A_{j,\varepsilon}^{-1} R_j.
\end{equation}
We define the matrices $A_{j,\varepsilon}$ in \eqref{bonazzoli_mini_10_eq:ORAS} to be
the matrices stemming from the discretization of the following local Robin problems with absorption 
\begin{align*}
	-\Delta u_j -(k^2+\mathtt{i}\varepsilon) u_j&= f &&\mbox{in } \Omega_j, \\
	\frac{\partial u_j}{\partial n_j} - \mathtt{i} \eta u_j &= 0 && \mbox{on } \partial \Omega_j.
\end{align*}

In order to achieve weak dependence on the wavenumber $k$ and number of subdomains, we add a coarse component to \eqref{bonazzoli_mini_10_eq:ORAS}. The \emph{two-level} preconditioner can be written in a generic way as follows 
\begin{equation}
	\label{bonazzoli_mini_10_eq:BNN}
	M^{-1}_{2,\varepsilon} = Q M^{-1}_{1,\varepsilon} P + Z E^{-1} Z^{*},
\end{equation}
where $*$ denotes the conjugate transpose, $M^{-1}_{1,\varepsilon}$ is the one-level
preconditioner \eqref{bonazzoli_mini_10_eq:ORAS}, $Z$ is a rectangular matrix with full column rank, $E=Z^{*}A_\varepsilon Z$ is the so-called coarse grid matrix, $\Xi=Z E^{-1} Z^{*}$ is the so-called coarse grid correction matrix. If $P=Q=I$ this is an \emph{additive} two-level preconditioner (which {would be} called two-level ImpRAS in \cite{graham:2015:domain}). 
If $P= I - A_\varepsilon \Xi$ and $Q= I - \Xi A_\varepsilon$, this is a \emph{hybrid} two-level preconditioner (ImpHRAS in \cite{graham:2015:domain}), also called the \emph{Balancing} Neumann Neumann (BNN) preconditioner. 
Preconditioner~\eqref{bonazzoli_mini_10_eq:BNN} is characterized by the choice of $Z$, whose columns span the \emph{coarse space} (CS). 
We will consider the following two cases:\\

\noindent {\bf The grid coarse space}
The most natural coarse space would be one based on a coarser mesh, we subsequently call it ``grid coarse space''. Let us consider $\mathcal T_{H_{\text{coarse}}}$  a simplicial mesh of $\Omega$ with mesh diameter $H_{\text{coarse}}$ and $V_{H_{\text{coarse}}} \subset V$ the corresponding finite
element space. 
\replaced[id=R2]{Let $\mathcal I_0\colon V_{H_{\text{coarse}}} \to V_h$ be the nodal interpolation operator and define $Z$ as the corresponding matrix. Then }{Let $\mathcal R_0\colon V_h \to V_{H_{\text{coarse}}}$ be the nodal interpolation operator 
and $R_0$ the corresponding matrix. Define $Z=R_0^T$, then} in this case $E=Z^{*}A_\varepsilon Z$ is really the stiffness matrix of the problem (with absorption) discretized on the coarse mesh. 
Related preconditioners without absorption are used in \cite{KimnSar:2007:ROB}.\\

\noindent {\bf The DtN coarse space}
This coarse space (see \cite{Conen:2014:CSH}) is based on local Dirichlet-to-Neumann (DtN) eigenproblems on the subdomain interfaces. 
For a subdomain $\Omega_i$, first of all consider $a^{(i)}\colon H^1(\Omega_i) \times H^1(\Omega_i) \to \mathbb{R}$
\begin{equation*}
	a^{(i)}(v,w) = \int_{\Omega_i}{\left(\nabla v \cdot \overline{\nabla w} - (k^2+\mathtt{i}\varepsilon) v \overline{w}\right)}
	-
	\int_{\partial \Omega_i \cap \partial \Omega}{\mathtt{i} \eta u \overline{v} }. 
\end{equation*}
Let $(A^{(i)})_{kl} = a^{(i)}\left(\phi_k, \phi_l\right)$, 
and let $I$ and $\Gamma_i$ be the sets of indices corresponding, resp., to the interior and
boundary dofs on $\Omega_i$, with $n_I$ and $n_{\Gamma_i}$ their cardinalities.
With the usual block notation, the subscripts $I$ and $\Gamma_i$ for the
matrices $A$ and $A^{(i)}$ denote the entries of these matrices associated with
the respective dofs.
Let
$
	M_{\Gamma_i} = \left( \int_{\Gamma_i}{ \phi_k \phi_l } \right)_{k,l \in \Gamma_i}
$
be the mass matrix on the interface $\Gamma_i=\partial \Omega_i \setminus \partial \Omega$ of subdomain $\Omega_i$.  We need to solve the following eigenproblem: find $\left(\vec u, \lambda\right) \in \mathbb
C^{n_{\Gamma_i}} \times \mathbb C$, s.t.
\begin{equation}
	\label{bonazzoli_mini_10_eq:DtNDiscEigenproblem}
	( A_{\Gamma_i \Gamma_i}^{(i)} - A_{\Gamma_i I} A_{II}^{-1} A_{I \Gamma_i}
	) \vec u = \lambda M_{\Gamma_i} \vec u. 
\end{equation}
Now, the matrix $Z$ of the DtN coarse space is a rectangular, block-diagonal matrix with blocks $W_i$, associated with the subdomain $\Omega_i$, $1\leq i \leq N_\text{sub}$, given by Algorithm~\ref{bonazzoli_mini_10_alg:WiConstruction}.
If $m_i$ is the number of eigenvectors selected by the automatic criterion in Line~\ref{bonazzoli_mini_10_alg:WiConstruction:Choice} of Algorithm~\ref{bonazzoli_mini_10_alg:WiConstruction}, the block $W_i$ has dimensions $n_i \times m_i$, and the matrix $Z$ has dimensions ${n \times \sum_{j=1}^{N_\text{sub}} m_i}$.
Due to the overlap in the decomposition, the blocks may share some rows inside the matrix $Z$.

\begin{algorithm}
  \caption{Construction of the block $W_i$  of the DtN CS matrix $Z$
    \label{bonazzoli_mini_10_alg:WiConstruction}}
  \begin{algorithmic}[1]
    \Statex
    	\State Solve the discrete DtN eigenproblem \eqref{bonazzoli_mini_10_eq:DtNDiscEigenproblem}
    	on subdomain $\Omega_i$ for the eigenpairs $(\lambda_j,\vec g_i^j)$.
	\label{bonazzoli_mini_10_alg:WiConstruction:DtN}
    	\State Choose the $m_i$ eigenvectors $\vec g_i^j \in
	\mathbb C ^{n_{\Gamma_i}}$ such that $\Re(\lambda_j)< k$, $1\leq j \leq m_i$. 
	\label{bonazzoli_mini_10_alg:WiConstruction:Choice}
		\For{$j = 1 \textrm{ to } m_i$}
			\State Compute the discrete Helmholtz extension
			$\vec u_i^j \in \mathbb C^{n_i}$ to $\Omega_i$ of $\vec g_i^j$ 
			as $\vec u_i^j  = [-A_{II}^{-1} A_{I \Gamma_i} \vec g_i^j, \quad \vec g_i^j]^T$.
			\label{bonazzoli_mini_10_alg:WiConstruction:Ext}
		\EndFor
		\State Define the matrix $W_i \in \mathbb
	C^{n_i \times m_i}$ as
	$W_i = 
		\begin{pmatrix}
			D_i \vec u_i^1, & \ldots, & D_i \vec u_i^{m_i}
		\end{pmatrix}$.
  \end{algorithmic}
\end{algorithm}

\section{Numerical experiments}

We solve \eqref{bonazzoli_mini_10_eq:matrixEqHelm} with $\varepsilon=0$ on the unit square/cube, with a uniform simplicial
mesh of diameter $h \sim k^{-3/2}$, which is believed to remove the pollution effect. 
The right-hand side is given by $f = -\exp({-100((x-0.5)^2+(y-0.5)^2)})$ for $d=2$, $f = -\exp(-400((x-0.5)^2+(y-0.5)^2+(z-0.5)^2))$ for $d=3$.

We use GMRES with right preconditioning (with a tolerance $\tau = 10^{-6}$), starting with a \emph{random initial guess}, 
\added[id=R1]{which ensures, unlike a zero initial guess, that all frequencies are present in the error;} 
the stopping criterion is based on the relative residual.
We consider a regular decomposition into subdomains (squares/cubes), the overlap for each subdomain is of size $\mathcal{O}(2h)$ in all directions and the two-level preconditioner~\eqref{bonazzoli_mini_10_eq:BNN} is used in the hybrid way.
All the computations are done in the open source language FreeFem++ (\url{http://www.freefem.org/ff++/}).
The $3d$ code is parallelized and run on the TGCC Curie supercomputer. 
We assign each subdomain to one processor.
So in our experiments the number of processors increases if the number of subdomains increases.
To apply the preconditioner, the local problems in each subdomain (with matrices $A_{j,\varepsilon}$ in \eqref{bonazzoli_mini_10_eq:ORAS}) and the coarse space problem (with matrix $E$ in \eqref{bonazzoli_mini_10_eq:BNN}) are solved with a direct solver. 

As in \cite{graham:2015:domain,Graham:2016:RRD}, in the experiments we take the subdomain diameter $H_{\text{sub}}$ and the coarse mesh diameter $H_{\text{coarse}}$ constrained by $k$: $H_{\text{sub}} \sim k^{-\alpha}$ and $H_{\text{coarse}} \sim k^{-\alpha'}$, for some choices of $0<\alpha,\alpha'<=1$ detailed in the following; if not differently specified, we take $\alpha=\alpha'$, which is the setting of all numerical experiments in \cite{graham:2015:domain}. Note that $H_{\text{coarse}}$ does not appear as a parameter in the DtN coarse space.
We denote by $n_\text{CS}$ the size of the coarse space. 
For the grid coarse space 
$n_\text{CS}=(1/H_{\text{coarse}}+1)^d$, the number of dofs for the nodal linear finite elements in the unit square/cube. For the DtN coarse space $n_\text{CS} = \sum_{j=1}^{N_\text{sub}} m_i$, the total number of  computed eigenvectors for all the subdomains.
While we solve the pure Helmholtz problem without absorption, both the one-level preconditioner~\eqref{bonazzoli_mini_10_eq:ORAS} and the two-level preconditioner~\eqref{bonazzoli_mini_10_eq:BNN} are built from problems which can have non zero absorption given by $\varepsilon_\text{prec}=k^\beta$. In the experiments we put $\beta=1$ or $\beta=2$.

In the following tables we compare the one-level preconditioner, 
the two-level preconditioners with the grid coarse space and with the DtN coarse space 
in terms of number of iterations of GMRES and size of the coarse space ($n_\text{CS}$), for different values of the wavenumber $k$ and of the parameters $\alpha, \beta$. 
We also report the number of subdomains $N_\text{sub}$, which is controlled by $k$ and $\alpha$ as mentioned above. Since $h \sim k^{-3/2}$, the dimension $n$ of the linear system matrix is of order $k^{3d/2}$; for 3d experiments we report $n$ explicitly.
Tables~\ref{bonazzoli_mini_10_tab:2dHelm_betas}, \ref{bonazzoli_mini_10_tab:2dHelm_sameCSsize} concern the $2d$ problem, Table~\ref{bonazzoli_mini_10_tab:3dHelm_alphaprime} the $3d$ problem.

\begin{table}
\begin{center}
\begin{tabular}{|c|c||c|c|c|c|c|}
\hline 
\multicolumn{2}{|c||}{} & \multicolumn{5}{|c|}{$\beta=1$}\tabularnewline
\hline
\multicolumn{2}{|c||}{} & \multicolumn{5}{|c|}{$\alpha=0.6$}\tabularnewline
\hline 
$k$ & $N_\text{sub}$ & 1-level & grid CS & $n_\text{CS}$ & DtN CS & $n_\text{CS}$\tabularnewline
\hline 
10 & 9 & 22 & 19 & 16 & 11 & 39\tabularnewline
20 & 36 & 48 & 46 & 49 & 26 & 204\tabularnewline
40 & 81 & 78 & 98 & 100 & 37 & 531\tabularnewline
60 & 121 & 109 & 114 & 144 & 43 & 1037\tabularnewline
80 & 169 & 139 & 138 & 196 & 93 & 1588\tabularnewline 
\hline
\multicolumn{2}{|c||}{} & \multicolumn{5}{|c|}{$\alpha=0.8$}\tabularnewline
\hline 
$k$ & $N_\text{sub}$ & 1-level & grid CS & $n_\text{CS}$ & DtN CS & $n_\text{CS}$\tabularnewline
\hline 
10 & 36 & 35 & 19 & 49 & 10 & 122\tabularnewline
20 & 100 & 71 & 35 & 121 & 13 & 394\tabularnewline
40 & 361 & 158 & 88 & 400 & 22 & 1440\tabularnewline
60 & 676 & 230 & 187 & 729 & 39 & 2700\tabularnewline
80 & 1089 & 304 & 331 & 1156 & 68 & 4352\tabularnewline
\hline
\multicolumn{2}{|c||}{} & \multicolumn{5}{|c|}{$\alpha=1$}\tabularnewline
\hline 
$k$ & $N_\text{sub}$ & 1-level & grid CS & $n_\text{CS}$ & DtN CS & $n_\text{CS}$\tabularnewline
\hline 
10 & 100 & 65 & 26 & 121 & 11 & 324\tabularnewline
20 & 400 & 122 & 26 & 441 & 14 & 1120\tabularnewline
40 & 1600 & 286 & 33 & 1681 & 20 & 4640\tabularnewline
60 & 3600 & 445 & 45 & 3721 & 29 & 10560\tabularnewline
80 & 6400 & $>$500 & 62 & 6561& 44 & 18880\tabularnewline
\hline 
\end{tabular}~
\begin{tabular}{|c|c|c|c|c|c|}
\hline 
\multicolumn{5}{|c|}{$\beta=2$}\tabularnewline
\hline
\multicolumn{5}{|c|}{ $\alpha=0.6$}\tabularnewline
\hline 
1-level & grid CS & $n_\text{CS}$ & DtN CS & $n_\text{CS}$\tabularnewline
\hline 
28 & 27 & 16 & 23 & 40\tabularnewline
67 & 56 & 49 & 40 & 220\tabularnewline
121 & 114 & 100 & 72 & 578 \tabularnewline
169 & 165 & 144 & 109 & 920\tabularnewline
223 & 216 & 196 &  126 & 1824 \tabularnewline 
\hline
\multicolumn{5}{|c|}{ $\alpha=0.8$}\tabularnewline
\hline 
1-level & grid CS & $n_\text{CS}$ & DtN CS & $n_\text{CS}$\tabularnewline
\hline 
39 & 27 & 49 & 28 & 86\tabularnewline
83 & 51 & 121 & 41 & 362\tabularnewline
182 & 95 & 400 & 71 &1370\tabularnewline
268 & 150 & 729 & 103 & 2698\tabularnewline
355 & 214 & 1156 & 138 & 4350\tabularnewline
\hline
\multicolumn{5}{|c|}{ $\alpha=1$}\tabularnewline
\hline 
1-level & grid CS & $n_\text{CS}$ & DtN CS & $n_\text{CS}$\tabularnewline
\hline 
57 & 30 & 121 & 23 & 324\tabularnewline
130 & 49 & 441 & 42 & 1120\tabularnewline
296 & 80 & 1681 & 72 & 4640\tabularnewline
455 & 112 & 3721 & 101 & 10560\tabularnewline
$>$500 & 149 & 6561 & 134 & 18880\tabularnewline
\hline 
\end{tabular}
\end{center}
\caption{($d=2$) Number of iterations (and coarse space size $n_\text{CS}$) for the one-level preconditioner and the two-level preconditioners with the grid coarse space/DtN coarse space, with  $H_{\text{sub}}=H_{\text{coarse}} \sim k^{-\alpha}$, $\varepsilon_\text{prec}=k^\beta$.}
\label{bonazzoli_mini_10_tab:2dHelm_betas} 
\end{table}
In Table~\ref{bonazzoli_mini_10_tab:2dHelm_betas}, we let the DtN coarse space size be determined by the automatic choice criterion in Line~\ref{bonazzoli_mini_10_alg:WiConstruction:Choice} of Algorithm~\ref{bonazzoli_mini_10_alg:WiConstruction} (studied in \cite{Conen:2014:CSH}) and the grid coarse space size by $H_{\text{coarse}} \sim k^{-\alpha}$. We see that the DtN coarse space is much larger than the grid coarse space and gives fewer iterations. 
The preconditioners with absorption $\varepsilon_\text{prec}=k^2$ perform much worse than those with absorption $\varepsilon_\text{prec}=k$ independently of $n_\text{CS}$.
For $\varepsilon_\text{prec}=k$, when $\alpha=1$ \replaced[id=R2]{the number of iterations grows as $k^{0.9}$, respectively $k^{1.1}$, for the grid coarse space, respectively DtN coarse space (excluding the first two values for $k$ small where the asymptotic behaviour is not reached yet)}{ the number of iterations grows mildly with the wavenumber $k$ for both coarse spaces (but at the cost of an increasing coarse space size),} while the one-level preconditioner performs poorly \added[id=R2]{(for $k=80$ it needs more than $500$ iterations to converge)}. When $\alpha<1$, i.e. for coarser coarse meshes, the growth with $k$ is higher, and for $\alpha=0.6$ the two-level preconditioner is not much better than the one-level preconditioner because the coarse grid problem is too coarse; for $\alpha=0.8$ with the DtN coarse space the growth with $k$ degrades less than with the grid coarse space.

\begin{table}
\begin{center}
\begin{tabular}{|c|c||c|c|c|c|c|}
\hline 
\multicolumn{2}{|c||}{} & \multicolumn{4}{c|}{$n_\text{CS}$ forced by grid CS}\tabularnewline
\hline
\multicolumn{2}{|c||}{} & \multicolumn{4}{c|}{$\alpha=0.6$}\tabularnewline
\hline 
$k$ & $N_\text{sub}$ & grid CS & $n_\text{CS}$ & DtN CS & $n_\text{CS}$\tabularnewline
\hline 
10 & 9 & 19 & 16 & 18 & 18\tabularnewline
20 & 36 & 46 & 49 &  44 & 72\tabularnewline
40 & 81 & 98 & 100 & 85 & 162\tabularnewline
60 & 121 & 114 & 144 & 109 & 242\tabularnewline
80 & 169 & 138 & 196 & 140 & 338\tabularnewline
\hline
\multicolumn{2}{|c||}{} & \multicolumn{4}{c|}{$\alpha=0.8$}\tabularnewline
\hline 
$k$ & $N_\text{sub}$ & grid CS & $n_\text{CS}$ & DtN CS & $n_\text{CS}$\tabularnewline
\hline 
10 & 36 & 19 & 49 & 26 & 72\tabularnewline
20 & 100 & 35 & 121 & 61 & 200\tabularnewline
40 & 361 & 88 & 400 & 139 & 722\tabularnewline
60 & 676 & 187 & 729 & 191 & 1352\tabularnewline
80 & 1089 & 331 & 1156 & 250 & 2178\tabularnewline
\hline
\multicolumn{2}{|c||}{} & \multicolumn{4}{c|}{$\alpha=1$}\tabularnewline
\hline 
$k$ & $N_\text{sub}$ & grid CS & $n_\text{CS}$ & DtN CS & $n_\text{CS}$\tabularnewline
\hline 
10 & 100 & 26 & 121 & 52 & 200\tabularnewline
20 & 400 & 26 & 441 & 43 & 800\tabularnewline
40 & 1600 & 33 & 1681 & 157 & 3200\tabularnewline
60 & 3600 & 45 & 3721 & 338 & 7200\tabularnewline
80 & 6400 & 62 & 6561& $>$500& 12800\tabularnewline
\hline 
\end{tabular}~
\begin{tabular}{|c|c|c|c|c|c|}
\hline 
\multicolumn{4}{|c|}{$n_\text{CS}$ forced by DtN CS}\tabularnewline
\hline
\multicolumn{4}{|c|}{ $\alpha=0.6$}\tabularnewline
\hline 
grid CS & $n_\text{CS}$ & DtN CS & $n_\text{CS}$\tabularnewline
\hline 
17 & 36 & 11 & 39\tabularnewline
24 & 196 & 26 & 204\tabularnewline
50 & 529 & 37 & 531\tabularnewline
104 & 841 & 43 & 1037 \tabularnewline
173 & 1521 & 93 & 1588 \tabularnewline
\hline
\multicolumn{4}{|c|}{ $\alpha=0.8$}\tabularnewline
\hline 
grid CS & $n_\text{CS}$ & DtN CS & $n_\text{CS}$\tabularnewline
\hline 
15 & 121 & 10 & 122\tabularnewline
20 & 361 & 13 & 394\tabularnewline
35 & 1369 & 22 & 1440\tabularnewline
52 & 2601 & 39 & 2700\tabularnewline
78 & 4225 & 68 & 4352 \tabularnewline
\hline
\multicolumn{4}{|c|}{ $\alpha=1$}\tabularnewline
\hline 
grid CS & $n_\text{CS}$ & DtN CS & $n_\text{CS}$\tabularnewline
\hline 
17 & 324 & 11 & 324\tabularnewline
23 & 1089 & 14 & 1120\tabularnewline
22 & 4624 & 20 & 4640\tabularnewline
26 & 10404 & 29 & 10560\tabularnewline
30 & 18769 & 44 & 18880 \tabularnewline
\hline 
\end{tabular}
\end{center}
\caption{($d=2$) Number of iterations (and coarse space size $n_\text{CS}$) for the two-level preconditioners with the grid coarse space/DtN coarse space \emph{forcing similar $n_\text{CS}$}, with  $H_{\text{sub}} \sim k^{-\alpha}$, $\varepsilon_\text{prec}=k$.}
\label{bonazzoli_mini_10_tab:2dHelm_sameCSsize} 
\end{table}
We have seen in Table~\ref{bonazzoli_mini_10_tab:2dHelm_betas} that the DtN coarse space gives fewer iterations than the grid coarse space, but their sizes differed significantly. Therefore, in Table~\ref{bonazzoli_mini_10_tab:2dHelm_sameCSsize} we compare the two methods forcing  $n_\text{CS}$ to be similar. On the left, we force the DtN coarse space to have a smaller size, similar to the one of the grid coarse space, by taking just $m_i=2$ eigenvectors for each subdomain. On the right, we do the opposite, we force the grid coarse space to have the size of the DtN coarse space obtained in~Table~\ref{bonazzoli_mini_10_tab:2dHelm_betas}, by prescribing a smaller coarse mesh diameter $H_{\text{coarse}}$, while keeping the same number of subdomains as in Table~\ref{bonazzoli_mini_10_tab:2dHelm_betas} with $H_{\text{sub}} \sim k^{-\alpha}$.
We can observe that for smaller coarse space sizes (left) the grid coarse space gives fewer iterations than the DtN coarse space, while for larger coarse space sizes (right) the result is reversed.

\begin{table}
\begin{center}
\begin{tabular}{|c|c|c||c|c|c|c|c|}
\hline
\multicolumn{3}{|c||}{} & \multicolumn{5}{|c|}{$\alpha=0.5$, $\alpha'=1$}\tabularnewline
\hline 
$k$ & $n$ & $N_\text{sub}$ & 1-level & grid CS & $n_\text{CS}$ & DtN CS & $n_\text{CS}$\tabularnewline
\hline 
10 & 39304 & 27 & 25 & 12 & 1331 & 14 & 316\tabularnewline
20 & 704969 & 64 & 39 & 17 & 9261 & 31 & 1240 \tabularnewline
30 & 5000211 & 125 & 55 & 21 & 29791 & 54 & 2482 \tabularnewline
40 & 16194277 & 216 & 74 & 29 & 68921 & 80 & 4318\tabularnewline
\hline
\multicolumn{3}{|c||}{} & \multicolumn{5}{|c|}{$\alpha=0.6$, $\alpha'=0.9$}\tabularnewline
\hline 
$k$ & $n$ & $N_\text{sub}$ & 1-level & grid CS & $n_\text{CS}$ & DtN CS & $n_\text{CS}$\tabularnewline
\hline 
10 & 39304 & 27 & 25 & 15 & 512 & 14 & 316\tabularnewline
20 & 912673 & 216 & 61 & 24 & 3375 & 41 & 2946\tabularnewline
30 & 4826809 & 343 & 73 & 34 & 10648 & 65 & 6226 \tabularnewline
40 & 16194277 & 729 & 98 & 48 & 21952 & 108 & 13653 \tabularnewline
\hline
\multicolumn{3}{|c||}{} & \multicolumn{5}{|c|}{$\alpha=0.7$, $\alpha'=0.8$}\tabularnewline
\hline 
$k$ & $n$ & $N_\text{sub}$ & 1-level & grid CS & $n_\text{CS}$ & DtN CS & $n_\text{CS}$\tabularnewline
\hline 
10 & 46656 & 125 & 34 & 19 & 343 & 11 & 896 \tabularnewline
20 & 912673 & 512 & 73 & 35 & 1331 & 18 & 4567 \tabularnewline
30 & 5929741 & 1000 & 103 & 57 & 4096 & 65 & 12756\tabularnewline
40 & 17779581 & 2197 & 139 & 89 & 8000 & 116 & 30603 \tabularnewline
\hline 
\multicolumn{3}{|c||}{} & \multicolumn{5}{|c|}{$\alpha=0.8$, $\alpha'=0.7$}\tabularnewline
\hline 
$k$ & $n$ & $N_\text{sub}$ & 1-level & grid CS & $n_\text{CS}$ & DtN CS & $n_\text{CS}$\tabularnewline
\hline 
10 & 50653 & 216 & 39 & 23 & 216 & 19 & 1354\tabularnewline
20 & 1030301 & 1000 & 46 & 86 & 729 & 23 & 7323\tabularnewline
30 & 5929741 & 3375 & 137 & 116 & 1331 & 21 & 26645\tabularnewline
40 & 28372625 & 6859 & 189 & 200 & 2744 & 27 & 54418\tabularnewline
\hline 
\end{tabular}
\end{center}
\caption{($d=3$) Number of iterations (and coarse space size $n_\text{CS}$) for the one-level preconditioner and the two-level preconditioners with the grid coarse space/DtN coarse space, with  $H_{\text{sub}}\sim k^{-\alpha}$, $H_{\text{coarse}} \sim k^{-\alpha'}$, $\varepsilon_\text{prec}=k$.}
\label{bonazzoli_mini_10_tab:3dHelm_alphaprime} 
\end{table}
We have seen that the coarse mesh obtained with $H_{\text{coarse}} \sim k^{-\alpha'}$, $\alpha'=\alpha$ can be too coarse if $\alpha=0.6$. At the same time, for $\alpha=1$ the number of subdomains gets quite large since $H_{\text{sub}} \sim k^{-\alpha}$, especially in $3d$; this is not desirable because in our parallel implementation we assign each subdomain to one processor, so communication among them would prevail and each processor would not be fully exploited since the subdomains would become very small. 
Therefore, to improve convergence with the grid coarse space while maintaining a reasonable number of subdomains, we consider separate $H_{\text{coarse}}$ and $H_{\text{sub}}$, taking $\alpha' \ne \alpha$.
For load balancing (meant as local problems having the same size as the grid coarse space problem), in $3d$ we choose $\alpha'=3/2-\alpha$. 
\replaced[id=R2]{The DtN coarse space  is still built by keeping the eigenvectors verifying the automatic choice criterion; note that in $3d$ the number of selected eigenvectors is larger than in $2d$, but we only keep a maximum of $20$ eigenvectors in each subdomain.}{The DtN coarse space size is still determined by the automatic choice criterion (among $20$ computed local eigenvectors) in each subdomain.} 
In Table~\ref{bonazzoli_mini_10_tab:3dHelm_alphaprime} we report the results of this experiment. 
As expected, for the grid coarse space the best iteration counts are obtained for $\alpha=0.5$ because then $\alpha'=1$ gives the coarse mesh with the smallest diameter among the experimented ones: the number of iterations grows slowly, with $\mathcal{O}(k^{0.61}) \cong \mathcal{O}(n^{0.13})$. With higher $\alpha$ the iteration counts get worse quickly, and $\alpha = 0.8$ is not usable.
For the DtN coarse space, the larger coarse space size is obtained by taking $\alpha$ bigger (recall that $\alpha'$ is not a parameter in the DtN case): for $\alpha=0.8$ the number of iterations grows slowly, with $\mathcal{O}(k^{0.2}) \cong \mathcal{O}(n^{0.04})$, but this value may be optimistic, there is a decrease in iteration number between $k = 20$ and $30$.
We believe that for the other values of $\alpha$, where the iteration counts are not much better or worse than with the one-level preconditioner, we did not compute enough eigenvectors in each subdomain to build the DtN coarse space. 




\section{Conclusion}

We tested numerically two different coarse space definitions for two-level domain decomposition preconditioners for the pure Helmholtz equation (discretized with piecewise linear finite elements), both in $2d$ and $3d$, reaching more than $15$ million degrees of freedom in the resulting linear systems.  
The preconditioners built with absorption $\varepsilon_\text{prec}=k^2$ appear to perform much worse than those with absorption $\varepsilon_\text{prec}=k$.
We have seen that in most cases for smaller coarse space sizes the grid coarse space gives fewer iterations than the DtN coarse space, while for larger coarse space sizes the grid coarse space gives generally more iterations than the DtN coarse space.
The best iteration counts for the grid coarse space are obtained by separating the coarse mesh diameter $H_{\text{coarse}} \sim k^{-\alpha'}$ from the subdomain diameter $H_{\text{sub}}\sim k^{-\alpha}$, taking $\alpha'>\alpha$. Both for \replaced[id=R2]{the grid coarse space}{ the coarse grid space} and the DtN coarse space, for appropriate choices of the method parameters we have obtained iteration counts which grow quite slowly with the wavenumber $k$.
Further experiments to compare \replaced[id=R2]{the two coarse spaces}{ the two definitions of coarse space} should be carried out in the heterogenous case.

\smallskip
\noindent \textbf{Acknowledgement} 
This work has been supported in part by the French National Research Agency (ANR), project MEDIMAX, ANR-13-MONU-0012.

\bibliographystyle{plainnat} 
\bibliography{bonazzoli_mini_10}



\end{document}